\documentclass[a4paper, 12pt]{article}

\usepackage{amssymb,amsbsy,amsmath,amsfonts,amssymb,amscd}
\usepackage{latexsym,eucal,exscale,epic,eepic,epsfig,xcolor}
\usepackage{mathrsfs}

\newtheorem{theorem}{Theorem}[section]
\newtheorem{lemma}[theorem]{Lemma}
\newtheorem{proposition}[theorem]{Proposition}
\newtheorem{corollary}[theorem]{Corollary}

\newcommand{\lra}{\longrightarrow}
\newcommand{\noi}{\noindent}
\newcommand{\Ap}{\mathscr A}

\newcommand{\RR}{\mathscr  R}

\newcommand{\Ch}{\mathscr C\! h}

\newcommand{\WAff}{W^{\rm Aff}}

\newcommand{\TT}{\mathscr T}

\title{Orbits of a symmetric subgroup in the Bruhat-Tits building}
\author{Paul  Broussous}

\date{\today}

\begin{document}
\maketitle

\begin{abstract}

 Let $G$ be a connected reductive algebraic group defined over a non-archimedean locally compact field $F$ of odd residue characteristic. Let $\theta$ be an $F$-rational involution of $G$ and $H$ be the reductive $F$-group $G^\theta$. We describe the orbits of the symmetric subgroup $H_F$ in the set of chambers of the reduced building of $G$ over $F$.

\end{abstract}

 \tableofcontents

 \bigskip

\section{Introduction} 

Let $F$ be a non-archimean locally compact field of residue caracteristic not $2$. Let $G$ be a connected reductive group defined over $F$ equipped with an $F$-rational involution $\theta$. We consider the reductive symmetric space $H_F \backslash G_F$, where $G_F$ is the group of $F$-rational points in $G$ and $H_F$ is the subgroup of $\theta$-fixed elements. 

 Let $K$ be a {\it good} maximal compact subgroup of $G_F$ (the fixator of a {\it special vertex} in the Bruhat-Tits building of $G_F$). Some versions of the Cartan decomposition (or polar decomposition) have been obtained for the symmetric space $H_F\backslash G_F$ (\cite{DS}, \cite{BO}), where by a Cartan decomposition one means a description of the $(H_F ,K)$-double cosets. Of course having such a decomposition is motivated by questions regarding harmonic analysis on the coset space $H_\backslash G_F$.
\medskip 

 This paper is motivated by the study of the Iwahori-spherical representations of $G_F$ contributing to harmonic analysis on $H_F \backslash G_F$ (the so called {\it distinguished} representations). Recall that these representations
 have been constructed by A. Borel (\cite{Bo}) using the geometry of the Bruhat-Tits building of $G_F$. If $I$ denotes the connected Iwahori subgroup of $G_F$, and $(\pi ,V)$ is an irreducible 
 Iwahori-spherical representation of $G_F$, so that the fixed vector space $V^I$ is non-trivial,  Borel gave a model of $\pi$ as a space of functions defined on the set of chambers in Bruhat-Tits building of $G_F$ with values in $V^I$.
 
  The case of the Steinberg representation was studied by the author and F. Court\`es (\cite{BC}, \cite{C}) and other Iwahori-Spherical representations where discussed in \cite{Br}. 
  \medskip
  
   So the aim of this paper is to describe the double coset set $H_F
\backslash G_F /I$. In fact we prefer taking a more geometric point of view by
instead considering the quotient $H_F \backslash \Ch_G$, where $\Ch_G$ denotes the $G_F$-set of chambers in the affine building of $G_F$ (note that we have a natural identification 
 $H_F \backslash G_F /I \simeq  H_F \backslash \Ch_G$, when $G$ is semisimple and simply connected).
 This quotient turns out to be much easier to understand than $H_F \backslash G_F /K$. 
 Its description is quite similar to that of $H_F \backslash G_F /P_{\rm min}$ where $P_{\rm min}$ is a minimal parabolic subgroup of $G_F$ (e.g. see \cite{HW} Proposition 6.10), and involves the set $H_F \backslash
 \mathscr T$ of $H_F$-conjugacy classes of maximal  $\theta$-stable $F$-tori in $G$ that split over $F$. 
 \medskip
 
  We prove the following result (see theorem \ref{decomposition} below) : 
  
  \begin{theorem} We have the disjoint union
  decomposition :
  	$$
H_F\backslash \Ch_G =
	\coprod_{T\in H_F \backslash {\mathscr T}} (N(T)_F\cap H_F)	\backslash \Ch_T \ .
	$$
\noi where  $\Ch_T$ denotes the set of chambers in  the apartment of a
torus $T$) and $N(T)_F$ the group of $F$-points of the normalizer of a torus $T$. 
 \end{theorem}
 
  We make this decomposition entirely explicit when $G$ is semisimple, simply connected and split over $F$ and when
  $\theta$ is a Galois involution attached  to a quadratic unramified extension (Proposition \ref{decompositionWeyl}). 
  For this we use a description of $H_F \backslash {\mathscr T}$ due to Paul G\'erardin (\cite{G}).
  \medskip
  
  The first ingredient of our proof is due to Delorme and S\'echerre (\cite{DS} Proposition 2.1) and was a key argument in loc. cit.: any chamber $C$ in $\Ch_G$, is contained in 
   a $\theta$-stable apartment  $\Ap$
   of the affine building of $G$. Our second ingredient 
   (Theorem \ref{unique}) 
   makes Delorme and S\'echerre's result more precise: the apartment $\Ap$
   is unique up to the action of an element of $P_{C,0+}\cap H_F$, where $P_{C,0+}$ is the pro-unipotent radical of the Iwahori subgroup $P_C$ fixing $C$. Note that this an analogue of the result stating that any minimal parabolic subgroup $P$ of $G_F$ contains a $\theta$-stable maximal split torus, which is unique up to the action of ${\rm Rad}(P)\cap H_F$, where ${\rm Rad}(P)$ is the unipotent radical of $P$ (\cite{HW} Lemma 2.4).
 
 \bigskip

{\bf N.B.} This morning Guy  Shtotland posted  the following article on MathArXiv: {\it Iwahori Matsumoto presentation for modules of Iwahori fixed functions on symmetric spaces} ( https://arxiv.org/abs/2406.16070 ). It contains a proof of my main theorem. Both articles, his and mine,  were written independently. 

\section{Notation for algebraic groups and buildings} 

 We let $F$ be a non-archimedean, non discrete,  locally compact field of {\bf odd} residue characteristic. 
If $M$ is an algebraic group defined over $F$, we denote by $M_F$ its group of $F$-rational points.

 We fix  once for all:
\medskip

 -- a reductive algebraic group $G$ that we assume connected and defined over $F$,

 --  an $F$-rational involution of $G$, that is a non trivial $F$-rational group automorphism $\theta$ of $G$ such that  $\theta^2 ={\rm id}$. 
\medskip

 We set $H=G^\theta$.  Thanks to the hypothesis on the residue characteristic of $F$, the group $H$ is a (non-necessarily connected) reductive group defined over $F$. 
\medskip

 We denote by $X={\rm BT}(G,F)$ the non-enlarged Bruhat-Tits building of $G$ over $F$. Recall that this is a metric space endowed with an action of $G_F$ by isometries. Moreover it has the structure of a polysimplicial complex on which $G_F$ acts by simplicial automorphism. Its maximal (open) polysimplices are called chambers and we denote by $\Ch_G$ the corresponding $G_F$-set.

 The automorphism $\theta$ induced a polysimplicial automorphism of $X$ that we shall denote by the same symbol $\theta$. This bijection permutes the chambers of $X$ and satisfies $\theta (g.x)=\theta (g).\theta (x)$, for all $g\in G$, $x\in X$. 
\medskip

 Attached to any maximal $F$-split torus $T$ of $G$ defined over $F$, we have the following objects:
\medskip

 -- Its (relative) root system $\Phi (G,T)$, whose affine Weyl group will be denoted by $W_T^{\rm Aff}$.

 -- Its centralizer $Z=Z(T)$ in $G$. Recall that it is a connected reductive group defined over $F$ such that $Z/T$ is anisotropic. In particular $Z_F$ contains a unique connected parahoric subgroup
  that we denote by $Z_F^0$. 

 -- Its normalizer $N=N(T)$

 
 -- Its apartment $\Ap =\Ap (T) \subset X$. It is acted upon by $N_F$ and $\WAff_T$. The parahoric subgroup $Z_F^0$ fixes $\Ap$ pointwise and the global stabilizer of $\Ap$ in $G_F$ is $N_F$.

 -- $\Ch_T$ the $N_F$-set of chambers contained in $\Ap$. Recall that $N_F$ acts transitively on $\Ch_T$. 
 \medskip
 
   Let $\Omega$ a subset of $X$ contained in some apartement $\Ap (T)$. We shall use the compact open subgroups $P_\Omega$ and $U_\Omega$ defined in \cite{BT}, {\S}7, p. 157 (these subgroups are respectively denoted by ${\mathscr P}_\Omega$ and $(G_F )^{\sharp}_{\Omega}$ in \cite{KP}, Definition 7.4.1, p. 255).  The group $P_\Omega$ is the connected {\it parahoric subgroup} attached to $\Omega$. It fixes $\Omega$ pointwise, but the pointwise fixator of $\Omega$ is $G_F$ is in general bigger than $P_\Omega$. Recall that we have the decomposition $P_\Omega =Z_F^0\, U_\Omega$, and that the subgroup $U_\Omega$ acts transitively on the set of apartments of $X$ containing $\Omega$ (\cite{KP}, Proposition 7.6.4, p. 265). 
 \medskip
 
  Let $x$ be a point of $X$. We denote by $(P_{x,r})_{r\in \RR_{\geqslant}}$ the {\it Moy-Prasad filtration} of the parahoric subgroup $P_x = P_{\{ x\}}$ (\cite{KP}, Definition 13.2.1, p. 468, where $P_{x,r}$ is denoted by ${\mathscr P}_{x,r}$). Following Moy and Prasad, we set $P_{x,0+} := \bigcup_{r>0} P_{x,r}$. According to \cite{KP}, Proposition 13.5.2, $P_{x,0+}$ is a pro-$p$-subgroup of $G_F$ (a projective limit of $p$-subgroups). It is generally called the {\it pro-unipotent radical} of $P_x$. 
  
   If $\mathscr F$ is a facet of $X$, $P_{x,0}$ and $P_{x,0+}$ do not depend on $x\in \mathscr F$; we denote these subgroups by $P_{{\mathscr F},0}$ and $P_{{\mathscr F},0+}$ respectively. 
 
When $x$ lies in a chamber, that is when $P_x$ is an {\it Iwahori subgroup} of $G_F$, $P_{x,0+}$ is a maximal pro-$p$-subgroup of $G_F^0$ (\cite{KP} Proposition 13.5.2, p. 476, cf.  \cite{KP} Definition 2.6.23 for the definition of the reductive subgroup $G^0 \subset G$). In particular, the pro-unipotent radical of $Z_F^0$, that we shall denote by $Z_F^{0+}$,  is a maximal pro-$p$-subgroup of $Z_F^0$. 

\section{Some conjugacy resuls}

  We denote by $\TT$ the $H_F$-set of maximal split $F$-tori of $G$ that are stable by $\theta$. 
If $T\in \TT$, then $Z(T)$, $N(T)$, $\Ap (T)$ and $\Ch_T$ are stable under the action of $\theta$. 

 The following key result is due to Delorme and S\'echerre (but was proved independently by Court\`es  (\cite{C} Proposition 4.1) for the case of a Galois involution, when $G$ is split over $F$). 

\begin{theorem} \cite{DS} Let $C$ be a chamber of $X$. Then there exists a $\theta$-stable maximal $F$-split torus $T\in \TT$ such that $C$ is contained in $\Ap (T)$.
\end{theorem}

 Note that such a torus $T$ is not unique. However the following result implies  that $T$ is uniquely determined modulo the conjugacy action of $H_F$ on $\TT$.

\begin{theorem} \label{unique}
 Let $C$ be a chamber in $X$ and $T_1$, $T_2$ be two $\theta$-stable maximal $F$-split tori of $G$ containing $C$. Then there exists an element $g\in P_{C,0+}\cap H_F$ $g.\Ap (T_1) = \Ap (T_2 )$. 
\end{theorem}

\noi {\it Remark}. This theorem is due to Court\`es when $G$ is split and $\theta$ is a Galois involution (\cite{C} Proposition 4.8). Our proof is inspired from his proof.
\medskip

\noi {\it Proof}. Let $\Omega$ be the {\it enclos} (cf. \cite{BT}, (7.1.2), p. 157) of $C\cup \theta (C)$. This is a $\theta$-stable convex subset of $\Ap (T_1)\cap \Ap (T_2 )$. Hence $U_\Omega$ is $\theta$-stable.  Since $U_\Omega$ acts transitively on the set of apartments containing $\Omega$,  there exists $g\in U_\Omega$ such that $g.\Ap (T_1 ) = \Ap (T_2 )$. Applying $\theta$ to this equality, we obtain $\theta (g) .\Ap (T_1) = \Ap (T_2 )$, whence $g^{-1}\theta (g) .\Ap (T_1 ) =\Ap (T_1 )$ and $g^{-1} \theta (g)\in U_\Omega \cap N(T_1 )_F$. Write $\Theta$ for the two element group $\{ 1 ,\theta\}$. Then $\sigma \mapsto g^{-1}\sigma (g)$ is a $1$-cocycle of $\Theta$ in $U_\Omega \cap N(T_1 )_F$.

\begin{lemma} Let $T$ be a maximal $F$-split torus of $G$ and $\Omega$ be a subset of $\Ap (T)$containing a chamber $C$. Then $U_\Omega\cap N(T)_F\subset Z_F^0$.
\end{lemma} 

\noi {\it Proof of the lemma}.  Let $x$ be any point of $C$. Since $C\subset \Omega$, we have $U_\Omega \subset U_C = U_x$. By \cite{KP} Definition 7.4.1(3) and Corollary 7.4.9, we have $U_\Omega\cap N(T)_F \subset U_x \cap N(T)_F = U_x \cap Z(T)_F$. Moreover, by \cite{KP} Lemma 7.4.3, we have that $U_x \cap Z(T)_F \subset Z(T)_F^0$, and we are done.

\medskip

 Since $\Omega$ contains a chamber, the subgroup $U_\Omega$ is a pro-$p$-subgroup. It follows from the lemma that $U_\Omega\cap N(T_1 )$ is a pro-$p$-subgroup of $Z_F^0$, whence is contained in the pro-unipotent radical  $Z_F^{0+}$ of $Z_F^0$.
\medskip

We need another result, due to Hakim and Murnaghan.

\begin{lemma} (\cite{HM} Proposition 2.12.) Let $M$ be a connected reductive $F$-group acted upon by a $F$-rational involution $\theta$. 
Let $K$ be a $\theta$-stable open subgroup of $M_F$ contained in the pro-unipotent radical $M_{x,0+}$, form some point $x$ of the Bruhat-Tits building of $M$. Then the non-abelian cohomology set $H^1 (\Theta ,K)$ is trivial, where $\Theta =\{ 1,\theta\}$. 
\end{lemma}

 By applying the lemma to $M=Z$ and $K=U_\Omega\cap N(T_1 )_F$, we have that the cocycle $\sigma \mapsto g^{-1}\sigma (g)$ is split. So there exists $z\in Z_F^{0+}$ such that $g^{-1}\theta (g) = z^{-1}\theta (z)$, that is $gz^{-1}\in H_F$. So replacing $g$ by $gz^{-1}\in U_\Omega Z_F^{0+}$, we can arrange $g\in H_F$, $g.C = C$ and $g.\Ap (T_1 )=\Ap (T_2 )$. Since $\Omega\subset C$, we hace $U_\Omega \subset U_C$, so we may choose $g\in U_C Z^{0+}_F \cap H$. Finally since $U_\Omega Z_F^{0+}$ is a pro-$p$-subgroup of $P_C$, it is contained in $P_{C,0+}$; the theorem follows. 
 
 \begin{corollary} Let $T$ be a $\theta$-stable maximal $F$-split torus of $G$ and $C_1$, $C_2$ chambers in $\Ap (T)$. Then if $C_2 \in H_F . C_1$, there exists $h\in N(T)_F \cap H$ such that $C_2 =h.C_1$. 
 \end{corollary}

\noi {\it Proof}.  Since there exists $h\in H_F$ such that $C_2 =h.C_1$, we have $C_1 \in \Ap (T)\cap h^{-1}. \Ap (T)$. By the previous theorem, there exists $k\in H_F \cap P_C$ such that $h^{-1} .\Ap (T)= k.\Ap (T)$. It follows that $hk\in N(T)_F \cap H_F$ and we are done by replacing $h$ by $hk$. 
\medskip

For $T\in {\mathscr T}$, write $H_F\backslash{\mathscr T}$ for a set of representatives for the action of $H_F$ on $\mathscr T$. For each $T\in {\mathscr T}/H_F$, let $(N(T)_F\cap H_F )\backslash \Ch_T$
be  a set of representatives for the action of $N(T)_F\cap H_F$ on $\Ch_T$.  We have proved the following result. 

\begin{theorem} \label{decomposition} A set of representatives for the action of $H_F$ on $\Ch_G$ is given by the disjoint union
	$$
H_F\backslash \Ch_G =
	\coprod_{T\in H_F \backslash {\mathscr T}} (N(T)_F\cap H_F)	\backslash \Ch_T \ .
	$$
\end{theorem}

 \section{The case of a Galois symmetric space}

 Let us assume that
  $G$ is semisimple and simply connected. Then for each $T\in \mathscr T$, $N(T)_F$ acts on $\Ap (T)$ through the quotient $N(T)_F /Z(F)_F^0$. This quotient is isomorphic to the affine Weyl group $W_T^{\rm Aff}$  of the affine root system $\Phi^{\rm Aff} (G,T)$ of $G$ relative to $T$.  The group $N(T)_F\cap H_F$ acts on $\Ap (T)$
  through the quotient  $(N(T)_F\cap H_F )/(Z_F^0 \cap H_F)$.
 Since $\Ap (T)$ is a principal homogeneous space under the action  of $W_T^{\rm Aff}$, the decomposition of Theorem \ref{decomposition} gives a bijection :
 
$$
H_F\backslash \Ch_G  \simeq
\coprod_{T\in H_F \backslash {\mathscr T}} \left[ (N(T)_F\cap H_F )/(Z_F^0 \cap H_F) \right]  \backslash W_T^{\rm Aff}\ .
$$

Now assume that the algebraic group $H$ is a semisimple, simply connected, defined and split over $F$. Let $E/F$ be an unramified quadratic field extension  and set $G={\rm Res}_{E/F} H$ (restriction of scalars). Finally let $\theta$ be the action on $G$ of the non trivial element of ${\rm Gal}(E/F)$. 

Then if $T\in {\mathscr T}$, we have $Z=T$ and it is classical, since $Z_F^0$ is a connected parahoric subgroup of $Z_F$,  that $H^1 ({\rm Gal}(E/F), Z_F^0)$ is trivial (cf. \cite{KP} Lemma 8.1.4).  Hence the natural map $(N(T)_F\cap H_F )/(Z_F^0 \cap H_F) \lra \left( W_T^{\rm Aff}\right)^{{\rm Gal}(E/F)}$ is a group isomorphism.
 Let $\theta_T$ be the involutive affine transformation of $\Ap (T)$ induced by $\theta$. Then $\left( W_T^{\rm Aff}\right)^{{\rm Gal}(E/F)}$ coincide with the centralizer $W_T^{{\rm Aff}, \theta_T}$ of $\theta_T$ in $W_T^{\rm Aff}$.
 
  We have proved the following.
  
  \begin{proposition} \label{decompositionWeyl}    With the assumptions and notation as above, we have the decomposition :
  	$$
  	H_F\backslash \Ch_G =
  	\coprod_{T\in H_F \backslash {\mathscr T}} W_T^{{\rm Aff}, \theta_T}\backslash \Ch_T \ .
  	$$
  	
  \end{proposition}
 
   In order to make this last decomposition more explicit, we now use the classification of the elements of $\mathscr T$ due to Paul G\'erardin (\cite{G}  {\S}III, Tores maximaux non ramifi\'es).
\medskip

 Fix an element $T_o\in {\mathscr T}$ that is split over $F$. Such tori exist and are conjugate under $H_F$.  We abreviate $N_{0,F} = N(T_0 )_F$. According to {\it loc. cit} Th\'eor\`eme 2, page 109, $H_F  \backslash {\mathscr T}$ is in natural bijection with the set of conjugacy classes of elements of order dividing  $2$ in $W_{T_o}^{\rm Aff}$. Let us explain how to associate to an  element $w\in W_{T_o}^{\rm Aff}$ of order dividing $2$  a torus $T_w \in {\mathscr T}$. 

 Let $\tilde w\in N_{0,F}$ be an element lifting $w$. Since $H$ is semisimple and simply connected, the Galois cohomology set $H^1 ({\rm Gal}(E/F), H(E) )$ is trivial and there exists $g\in H(E)=G_F$ such that $g^{-1}\theta (g)={\tilde w}$. Then we set $T_w  = gT_o g^{-1}$. The $H_F$-conjugacy class of $T_w$ depends only on that of $w$ in $W_{T_o}^{\rm Aff}$.  Of course if $w=1$, we have $T_w = T_o$.
   
 The group  isomorphism $N_{o,F}\lra N(T_w )_F$, $n\mapsto gng^{-1}$, induces a group isomorphism $W_{T_o}^{\rm Aff}\lra W_{T_w}^{\rm Aff}$. This isomorphism is compatible with the affine bijection: $\Ap_{T_o}\mapsto \Ap_{T_w}$, $x\mapsto g.x$.  Through these identifications, the action of $\theta_{T_w}$ on $\Ap_{T_w}$ corresponds to the action of $w$ on $\Ap_{T_0}$. Hence denoting by $W_{T_o}^{{\rm Aff}, w}$ the centralizer of $w$ in $W_{T_o}^{\rm Aff}$, we obtain an identification $W_{T_o}^{{\rm Aff}, w}\simeq W_{T_w}^{{\rm Aff}, \theta_{T_w}}$.  As a consequence the decomposition of Proposition \ref{decompositionWeyl} may be rewritten:
$$
H\backslash \Ch_G \simeq \coprod_{w} W_{T_o}^{{\rm Aff}, w} \backslash \Ap_{T_o}\simeq \coprod_{w} W_{T_o}^{{\rm Aff}, w} \backslash W_{T_o}^{\rm Aff}\ ,
$$
\noi where $w$ runs over the $W_{T_o}^{\rm Aff}$-conjugacy classes of elements of $W_{T_o}^{\rm Aff}$ of order dividing $2$.

Paul Broussous
\smallskip

paul.broussous{@}math.univ-poitiers.fr
\medskip

 Laboratoire de Math\'ematiques et Applications, UMR 7348 du CNRS

Site du Futuroscope -- T\'el\'eport 2
 
11,  Boulevard Marie et Pierre Curie

B\^atiment H3 --  TSA 61125

86073 POITIERS CEDEX 

 France

\end{document}